\documentclass[12pt,reqno]{amsart}
\usepackage[latin1]{inputenc}
\usepackage[T1]{fontenc}
\usepackage{mathptmx}
\usepackage{amssymb}
\usepackage{tikz}
\usetikzlibrary{calc}

\usepackage{color}

\newtheorem{teo}{Theorem}[section]
\newtheorem{lema}[teo]{Lemma}
\newtheorem{prop}[teo]{Proposition}
\newtheorem{cor}[teo]{Corollary}

\newtheorem{defi}[teo]{Definition}
\newtheorem{ejem}[teo]{Example}
\newtheorem{algo}[teo]{Algorithm}

\newtheorem{nota}[teo]{Note}
\newtheorem{rem}[teo]{Remark}
\numberwithin{equation}{section}

\newenvironment{dem}{\noindent{\it Proof}.}{\qed}

\newcommand{\ff}{\mathbb{F}}
\newcommand{\pp}{\mathbb{P}_{\ff}^{2}}

\newcommand{\NNA}{\textbf{N}^{\ast}}
\newcommand{\NN}{\textbf{N}}

\newcommand{\ad}{\mathcal{A}}
\newcommand{\oo}{\mathcal{O}}
\newcommand{\ooo}{\overline{\mathcal{O}}}

\newcommand{\wchi}{\widetilde{\chi}}

\newcommand{\pe}{\underline{P}}
\newcommand{\me}{\underline{m}}
\newcommand{\mep}{\underline{mP}}
\newcommand{\ei}{\varepsilon_i}
\newcommand{\eu}{\varepsilon_1}
\newcommand{\ed}{\varepsilon_2}
\newcommand{\er}{\varepsilon_{r-1}}
\newcommand{\wss}{\Gamma_{\pe}}

\newcommand{\mepu}{(\underline{m}-\underline{1}) \underline{P}}
\newcommand{\uno}{\underline{1}}

\usepackage{pstricks}

\definecolor{light}{gray}{0.9}
\definecolor{medium}{gray}{0.8}

\let\epsilon=\varepsilon
\let\tilde=\widetilde

\begin{document}
\title{On multi-index filtrations associated to Weierstra\ss~semigroups}

\author{Julio Jos\'e Moyano-Fern\'andez}

\address{Universit\"at Osnabr\"uck, FB Mathematik/Informatik, 49069
Osnabr\"uck, Germany} \email{jmoyano@uos.de}

\thanks{The author was partially supported by the Spanish Government Ministerio de Educaci\'on y Ciencia (MEC), grants MTM2007-64704 and MTM2012--36917--C03--03 in cooperation with the European Union in the framework of the founds ``FEDER''}

\begin{abstract}
The aim of this paper is to review the main techniques in the computation of Weierstra\ss~semigroup at several points of curves defined over perfect fields, with special emphasis on the case of two points. Some hints about the usage of some packages of the computer algebra software \textsc{Singular} are also given.
\end{abstract}

\maketitle

\section{Introduction}

There are several classical problems in the theory of algebraic curves which are interesting from a computational point of view. One of them is the computation of the Weierstra\ss~semigroup of a smooth projective algebraic curve $\tilde{\chi}$ at a certain rational point $P$, together with a rational function $f_m \in \mathbb{F}(\tilde{\chi})$ regular outside $P$ and achieving a pole at $P$ of order $m$, for each $m$ in this semigroup. This problem is solved with the aid of the adjunction theory for plane curves, profusely developed by A. von Brill and M. Noether in the 19th century (see \cite{brino}, \cite{noe}) so that we assume the knowledge of a singular plane birational model $\chi$ for the smooth curve $\tilde{\chi}$.
\medskip

Given a smooth projective algebraic curve $\chi$ (over a perfect
field $\ff$) and a set $P_1, \ldots , P_r$ of (rational) points of
$\chi$, we consider the family of finitely dimensional vector
subspaces of $\ff (\chi)$ given by
$\mathcal{L}(\mep)=\mathcal{L}(m_1P_1+m_2P_2+ \ldots + m_rP_r)$,
where $\underline{m} = (m_1,\ldots, m_r) \in \mathbb{Z}^r$. This
family gives rise to a $\mathbb{Z}^r$-multi-index filtration on the
$\ff$-algebra $A$ of the affine curve $\chi \setminus \{P_1,
\ldots, P_r\}$, since one has $A=\bigcup_{\underline{m} \in
\mathbb{Z}^r} \mathcal{L}(\mep)$. This multifiltration is related
to Weierstra\ss~semigroups (with respect to several points in
general, see Delgado \cite{de}) and, in case of finite fields, to the
methodology for trying to improve the Goppa estimation of the
minimal distance of algebraic-geometrical codes, see for instance Carvalho and Torres \cite{cato}.
A connection of that
filtration with global geometrical-topological aspects in a
particular case was shown by Campillo, Delgado and Gusein-Zade \cite{cadegz}. Poincar\'e series associated to this filtrations in particular cases were studied by the author in \cite{moyano}
\medskip

Thus, a natural question is to provide a computational method in order to
estimate the values of $\dim_{\ff} \mathcal{L}(\mep) = \ell
(\mep)$ for $\underline{m} \in \mathbb{Z}^r$. More precisely, it
would be convenient to estimate and compute values of type $\ell
((\underline{m}+\underline{\varepsilon})\underline{P})-\ell(\mep)$
where $\underline{\varepsilon} \in \mathbb{Z}^r$ is a vector whose
components are $0$ or $1$. This can be done by extending the method
developed by Campillo and Farr\'an \cite{cafa} in the case $r=1$, based on the
knowledge of a plane model $\widetilde{\chi}$ for $\chi$ (with
singularities) and representing the global regular differentials
in terms of adjoint curves to $\widetilde{\chi}$. 
\medskip

The paper is organised as follows: Sections 2 and 3 are devoted to fix the algebraic-geometrical
prerequisites. Section 4 deals with the study of 
more specific questions concerning to our
purpose, namely the adjunction theory of curves, with the remarkable Brill-Noether
Theorem.
In Section 5 we define the Weierstra\ss~semigroup at several points and describe two methods to compute values
of the form $\ell
((\underline{m}+\underline{\varepsilon})\underline{P})-\ell(\mep)$.
The last section is devoted to
show and explain some procedures implemented in
\textsc{Singular} based on Section 5.
\medskip

Notice the practical relevance of these ideas in view of the
algebraic-geometric coding theory: the
Weierstra\ss~semigroup plays an important role in the decoding
procedure of Feng and Rao, see e.g. Campillo and Farr\'an \cite{cafa2}, or H\o holdt, van Lint and Pellikaan \cite{holipe}.

\section{Terminology and notations}

Let $\mathbb{F}$ be a perfect field, and let $\overline{\mathbb{F}}$ a fixed algebraic closure of $\ff$. 
Let $\chi$ be an absolutely irreducible projective algebraic curve
defined over $\ff$. We distinguish three types of points on $\chi$, namely the geometric points, i.e. those with coordinates on $\overline{\ff}$; the rational points, i.e. those with coordinates on $\ff$; and the closed points, which are residue classes of geometric points under the action of the Galois group of the field extension $\overline{\ff}/\ff$, namely
\[
P:=\{\sigma (p) : \sigma \in \mathrm{Gal} (\overline{\ff}/\ff)\},
\]
where $p$ is a geometric point. Notice that closed points correspond one to one to points on the curve $\chi$ viewed as an $\ff$-scheme
which are closed  for the Zariski topology. Every closed point has an associated residue field $\ff '$ which is a finite extension of $\ff$. The degree of a closed point $P$ is defined as the cardinal of its conjugation class, which equals the degree of the extension $\ff ' / \ff$. In particular, $P$ is rational if and only if $\deg P=1$.
\medskip

Let us assume $\chi$ to be non-singular (or, equivalently, smooth, since $\ff$ is perfect). Let $\overline{\mathbb{F}}(\chi)$ be the field of rational functions of $\chi$. Let $P$ be a closed point on $\chi$. The local ring $\mathcal{O}_{\chi,P}$ of $\chi$ at $P$ with maximal ideal $\mathfrak{m}_{\chi, P}$ is therefore a discrete valuation ring with associated discrete valuation $v_P$. An element $f \in \mathcal{O}_{\chi,P}$ is said to vanish at $P$ (or to have a zero at $P$) if $f \in \mathfrak{m}_{\chi, P}$. A rational function $f$ such that $f \notin \mathcal{O}_{\chi,P}$ is said to have a pole at $P$. The order of the pole of $f$ at $P$ is given by $|v_P(f)|$.
\medskip

A \emph{rational divisor} $D$ over $\ff$ is a finite linear
combination of closed points $P \in \chi$ with integer
coefficients $n_P$, that is, $D=\sum_P n_P P$. If $n_P \ge 0$ for
all $P$, then $D$ is called \emph{effective}. We define the
\emph{degree} of $D$ as $\deg D := \sum_P n_P \deg P$, and the
\emph{support} of $D$ as the set $\mathrm{supp}(D)=\big \{P \in
\chi \ \mathrm{closed} \mid n_P \ne 0  \big \}$.
The set of rational divisors on $\ff$ form an abelian group
$\mathcal{D}(\ff)$. Rational functions define \emph{principal} divisors, namely divisors of the form
\[
(f):=\sum_P \mathrm{ord}_P (f) P.
\]

A rational divisor $D=\sum n_P P$ defines a
$\ff$-vector space
\[
\mathcal{L}(D)=\Big \{f \in \ff(\chi)^{\ast} \mid (f) \ge -D \Big
\} \cup \big \{ 0 \big \},
\]
that is, the set of rational functions $f$ with poles only at the
points $P$ with $n_P \ge 0$ (and, furthermore, with the pole order
of $f$ at $P$ must be less or equal than $n_P$), and if $n_P<0$
such functions must have a zero at $P$ of order greater or equal
than $n_P$. The dimension $\ell (D):=\dim_{\ff} \mathcal{L}(D)$ is finite. 
Two elements $f,g \in \mathcal{L}(D)$ satisfy $(f)+D = (g)+D$ if
and only if $f=\lambda g$, $\lambda \in \ff$, i.e., if and only if
$f = \lambda g$ for a constant $\lambda \in \ff$. Therefore the set $|D|$ of effective
divisors equivalent to $D$ can be identified with the projective space
$\mathbb{P}_{\mathcal{L}(D)}$ of dimension $\ell (D)-1$. The set
$|D|$ is called a \emph{complete linear system} of $D$.
\medskip

Let $\Omega_{\ff}(\ff(\chi))$ be the module of differentials on $\ff(\chi)$. A differential form $\omega \in \Omega_{\ff}(\ff(\chi))$
defines a divisor $(\omega):=\sum_P \mathrm{ord}_P(\omega)P$, called a canonical divisor. A rational divisor $D$ defines again a $\ff$-vector space
\[
\Omega (D):= \{\omega \in \Omega_{\ff}(\ff(\chi))^{\ast} \mid
(\omega) \ge D \} \cup \{ 0 \}.
\] 
of finite dimension, denoted by $i(D)$. It is a central result in the theory of algebraic curves the interplay of the dimensions $\ell (D)$ and $i(D)$. The dimension $\ell (D)$ is bounded in the following sense:

\begin{prop}[Riemann's inequality]
There exists a nonnegative integer $g$ such that
\[
\ell (D) \ge \deg D +1 -g.
\]
for any
rational divisor $D$ on $\chi$.
\end{prop}

\begin{defi}
The smallest integer $g$ satisfying the Riemann's inequality
is called the \emph{genus} of $\chi$.
\end{defi}

Riemann's inequality tells us that if $D$ is a large divisor,
$\mathcal{L}(D)$ is also large. But we can be a bit more
precise by using $i(D)$:

\begin{teo}[Riemann-Roch]
Let $D$ be a rational divisor. Then:
$$
\ell (D) - i(D) = \deg D +1 -g.
$$
\end{teo}

\section{Rational parametrizations}

Let $\ff$ be a perfect field, and let $\chi$ be an absolutely irreducible algebraic plane curve defined over $\ff$. Let $P$ be a closed point on $\chi$.
Let us consider the local ring $\oo :=
\oo_{\chi,P}$ with maximal ideal $\mathfrak{m}$, and write $\ooo$ for the semilocal ring of the normalisation of $\chi$ at $P$. Finally, let 
$\hat{\oo}$ be the completion of $\oo$ with respect to the $\mathfrak{m}$-adic topology. Each maximal ideal of $\ooo$ (or, equivalently, every minimal prime ideal $\mathfrak{p}$ of $\hat{oo}$) is said to be a \emph{branch} of $\chi$ at $P$.
\medskip

Let us now choose an affine chart containing $P$ so that the curve $\chi$ has an equation $f(X,Y)=0$, and set $A:=\ff[X,Y]/(f(X,Y))$ as the affine coordinate ring. Notice that $\oo = A_P$. Hence
\[
\ff \subseteq \ff[X,Y]/(f(X,Y))=A \subseteq A_P=\oo.
\]
Since $\ff$ is perfect, we can apply Hensel's lemma to find a finite field extension $K / \ff$ such that $K \subseteq \hat{A_P}=\hat{\oo}$ is a coefficient field for $\hat{\oo}$. Moreover, $K$ is the integral closure of $\ff$ in $\hat{\oo}$.
\medskip

Since $\hat{\oo}\subseteq \overline{\hat{\oo}}\cong \hat{\ooo}$, one has
\[
K \subseteq \hat{\oo}/\mathfrak{p} \subseteq \overline{ \hat{\oo}/\mathfrak{p}} = \hat{\ooo_{\mathfrak{m}}},
\]
and we can apply Hensel's lemma again to obtain a finite extension $K' / K$ which is a coefficient field for the local ring $\hat{\ooo_{\mathfrak{m}}}$. 
Without loss of generality we can consider $P$ as the ideal $(X,Y)$ in $K[\![X,Y]\!]$ so that
$\hat{\oo}\cong K[\![X,Y]\!]/(f(X,Y))$. This implies the existence of natural morphisms
\[
K[\![X,Y]\!]/(f(X,Y)) \cong \hat{\oo} \longrightarrow \hat{\oo}/\mathfrak{p} \longrightarrow K'[\![t]\!] \cong \hat{\ooo_{\mathfrak{m}}}
\]
for any local uniformizing parameter $t \in \mathfrak{m} \setminus \mathfrak{m}^2$. Notice that $K$ can be considered isomorphic to the residue field at $P$. Preserving these notations, a \emph{parametrization} of the curve $\chi$ at the point $P$ related to the coordinates $X,Y$ is a $K$-algebra morphism $\rho: K[\![X,Y]\!] \longrightarrow K'[\![t]\!]$ being continuous for the $(X,Y)$-adic  and $t$-adic topologies and satisfying $\mathrm{Im}(\rho) \not\subseteq K'$ and $\rho (f)=0$. This is equivalent to give formal power series $x(t),y(t) \in K'[\![t]\!]$ with $x(t) \neq 0$ or $y(t)\neq 0$ such that $f(x(t),y(t))\equiv 0$. 
\medskip

Consider parametrizations $\rho:K[\![X,Y]\!] \to K'[\![t]\!]$ and $\sigma: K[\![X,Y]\!] \to K''[\![t]\!]$ of the same rational branch. The parametrization $\sigma$ is said to be \emph{derivated} from $\rho$ if there is a formal power series $\tau(u) \in K''[\![u]\!]$ with positive order and a continuous $K$-algebra morphism $\alpha: K'[\![t]\!] \to K''[\![u]\!]$ with $\alpha(t)=\tau(u)$ such that $\sigma=\alpha \circ \rho$. We write $\sigma \succ \rho$. The relation $\succ$ is a partial preorder. Two parametrizations $\sigma$ and $\rho$ are called \emph{equivalent} if $\sigma \succ \rho$ and $\rho \succ \sigma$.  Those parametrizations being minimal with respect to $\succ$ up to equivalence are called \emph{primitive}. Equivalent primitive parametrizations are called \emph{rational}. They always exist and are invariant under the action of the Galois group of the extension $\overline{K}/K$. Rational parametrizations are in one to one correspondence with rational branches of the curve (cf. Campillo and Castellanos \cite{caca}).

\section{Brill-Noether theory for curves}

This section contains a summary of the classic adjunction
theory of curves, started by Riemann \cite{rie} and developed by
M. Noether and A. von Brill in the 19th century. 
\medskip

Let $P$ be a closed point. Let $\mathcal{C}_P$ be the annihilator of the
$\oo$-module $\overline{\oo} / \oo$, i.e.
\[
\mathcal{C}_P=\mathcal{C}_{\overline{\oo}/\oo} = \{\varphi \in
\overline{\oo} \mid \varphi \overline{\oo} \subseteq \oo  \}.
\]
This set is the largest ideal in $\oo$ which is also an ideal in
$\overline{\oo}$, and it is called the \emph{conductor
ideal} of the extension $\overline{\oo} / \oo$. Since $\ooo$ is a semilocal Dedekind domain with maximal ideals 
$\overline{\mathfrak{m}}_{Q_1}, \ldots , \overline{\mathfrak{m}}_{Q_d}$ (where $Q_i$ denote the rational branches of $\chi$ at $P$), the conductor ideal has a unique factorisation
\[
\mathcal{C}_P=\prod_{i=1}^d \overline{\mathfrak{m}}_{Q_i}^{d_{Q_i}}
\]
as ideal in $\ooo$. The exponents $d_{Q_i}$ can be easily computed by means of the Dedekind formula (see Zariski \cite{za}): if $(x_i(t_i),y_i(t_i))$ is a rational parametrisation of $Q_i$ one has
\begin{equation} \label{eqn:dedekind}
d_{Q_i} =\mathrm{ord}_{t_{Q_i}}
\Bigg(\frac{f_Y(X(t_{Q_i}),Y(t_{Q_i}))}{X^{\prime}(t_{Q_i})}
\Bigg)= \mathrm{ord}_{t_{Q_i}}
\Bigg(\frac{f_X(X(t_{Q_i}),Y(t_{Q_i}))}{Y^{\prime}(t_{Q_i})}
\Bigg).
\end{equation}

Let $n:\wchi \to \chi$ ne the normalisation morphism of $\chi$. Notice that $\wchi$ is nonsingular with $\ff(\wchi)=\ff(\chi)$.  Let $\oo=\oo_{\chi, P}$ and $\ooo$ its normalisation. Let $Q \in n^{-1}(\{P\})$. Since $Q$ is nonsingular, it is $\mathcal{C}_P \cdot \oo = \mathfrak{m}_Q^{d_Q}$ for a nonnegative integer $d_Q$. We define the effective divisor
\[
\ad := \sum_P \sum_{Q \in n^{-1}(\{P\})} d_Q \cdot Q
\]
which is called the \emph{adjunction divisor} of $\chi$. Notice that $\ad$ is a well-defined divisor on $\wchi$ (in fact, if $P$ is nonsingular, there is only one $Q \in n^{-1}(\{P\})$  and in this case $d_Q=0$). This implies in particular that  the support of $\ad$ consists of all rational branches of $\chi$ at singular points. Moreover, by setting $n_P:=\dim_{\ff}\ooo/\mathcal{C}_P$ we have
\[
n_P=\sum_{Q \in n^{-1}(\{P\})} d_Q
\]
for every $P$ on $\chi$. Therefore $\deg \ad = \sum_{P \in \chi} n_P$ (cf. Arbarello et al. \cite[Appendix A]{arb}; also Tsfasman and Vl\u adu\c t \cite[2.5.2]{tsfas}).
\medskip

Let $F:=F(X_0,X_1,X_2)$ be a homogeneous (absolutely irreducible) polynomial of degree $d$ over $\ff$ which defines the projective plane curve $\chi$. Let $\mathcal{F}_d$ be the set of all homogeneous polynomials in three variables of degree $d$. Let $i:\chi \to \mathbb{P}^2_{\ff}$ be the embedding of $\chi$ into the projective plane and $\NN: \widetilde{\chi} \to
\mathbb{P}_{\ff}^{2}$ be the natural morphism given by $\NN = i
\circ n$. A rational divisor $D$ on $\mathbb{P}^2_{\ff}$ such that $\chi$ is not contained in $\mathrm{supp}(D)$ is called an \emph{adjoint divisor} of $\chi$ if the pull-back divisor $\NNA D$ satisfies $\mathrm{supp}(\ad) \subseteq \mathrm{supp}(\NNA D)$ for $\ad$ the adjunction divisor of $\chi$. We can consider the analogous notion at the level of homogeneous polynomials. For $H \in \mathcal{F}_d$ with $F \nmid H$ one can consider the pull-back $\NNA H$, which is actually the intersection divisor on $\wchi$ cut out by the plane curve defined by $H$ on $\mathbb{P}^2_{\ff}$, namely
\begin{equation}\label{eq:ad1}
\NNA H = \sum_{Q \in \wchi} r_Q \cdot Q,
\end{equation}
with $r_Q=\mathrm{ord}_Q (h)$ for $h \in \oo_{\chi, n(Q)}$ being a local equation of the curve defined by $H$ at the point $n(Q)$. If $H$ satisfies additionally $\NNA D \geq \ad$, then it will be called an \emph{adjoint form} on $\chi$, and the curve defined by $H$ will be called an \emph{adjoint curve} to $\chi$. Notice that adjoint curves there always exists (take for instance the polars of the curve, cf. Brieskorn and Kn\"orrer \cite{brie}, p. 599). 
\medskip

Let $d:=\deg \chi$. The differentials gob ally defined at $\chi$ are in one to one correspondence with adjoint curves on $\wchi$ of degree $d-3$:

\begin{teo}\label{teo:1}
Let  $\ad_n$ be the set of adjoints of degree $n$ of the curve
$\chi$ embedded in $\pp$, let $K_{\widetilde{\chi}}$ be a canonical
divisor on $\widetilde{\chi}$ and set $d := \deg \chi$.
For $n=d-3$ there is an $\ff$-isomorphism of complete linear
systems
\[
\begin{array}{ccc}
\ad_n& \longrightarrow& |K_{\widetilde{\chi}}|\\
D& \longmapsto& \NNA D - \ad.
\end{array}
\]
\end{teo}

The key idea is to realise that the map is injective since $n=d-3<d$; see Gorenstein \cite[p. 433]{gor} or \cite[2.2.1]{tsfas} for further details.
 \medskip
 
In practice, we know a priori the equation of the plane curve
$\chi$ (defined over a perfect field $\ff$) given by the form $F
\in \mathcal{F}_d$ and the data of a certain divisor
$R=\sum_{Q^{\prime}}r_{Q^{\prime}} \cdot Q^{\prime}$ (for finitely
many points $Q^{\prime}$ on $\wchi$) which is effective and
rational over $\ff$, involving a finite number of rational
branches $Q$ of $\chi$ and their corresponding coefficients.
Moreover, we are able to compute the adjunction divisor of $\chi$,
$\ad=\sum_Q d_Q \cdot Q$. Our aim is to interprete the
condition of being an adjoint form---called \emph{adjoint
condition}---given by (\ref{eq:ad1}) in terms of equations. More
generally, we are interesting in finding some \emph{adjoint} form
$H \in \ff[X_0,X_1,X_2]$ satisfying
\begin{equation} \label{eq:ad2}
\NNA H \ge \ad + R.
\end{equation}
This process is known as \emph{computing adjoint forms with base
conditions} (see \cite{cafa}, \S 4).\vspace{0.25cm}

First of all, we choose a positive integer $\widetilde{n} \in \mathbb{N}$ in such a way that
there exists an adjoint of degree $\widetilde{n}$ not being a
multiple of $F$ and satisfying
(\ref{eq:ad2}). A bound for $\widetilde{n}$ can be found in Hach\'e
\cite{hachetesis}. Take then also a form $H \in
\mathcal{F}_{\widetilde{n}}$ in a general way, what is nothing
else but taking a homogeneous polynomial in three variables of
degree $\widetilde{n}$ with its coefficients as indeterminates (that is,
$H(X_0,X_1,X_2)=\sum_{i+j+k=\widetilde{n}} \lambda_{i,j,k}
X_0^iX_1^jX_2^k$).  Second we compute a rational primitive parametrization $\big(
X(t),Y(t) \big)$ of $\chi$ at every branch involved in the
support of the adjunction divisor $\ad$ and the divisor $R$.
Next we get the support of the adjunction divisor $\ad$ from the
conductor ideal via the Dedekind formula (\ref{eqn:dedekind}). 
Last we consider the coefficient $r_Q$ of the
divisor $R$ at $Q$, and thus the local condition at $Q$ imposed on
$H$ by (\ref{eq:ad2}) is given by
\begin{equation} \label{eq:dqrq}
\mathrm{ord}_t h(X(t),Y(t)) \ge d_Q + r_Q,
\end{equation}
with $h$ the local affine equation of $H$ at $Q$. 
The inequality (\ref{eq:dqrq}) expresses a linear condition (given
by a linear inequation) on the coefficients $\lambda_{i,j,k}$ of
$h$. \vspace{0.25cm}

The required linear equations are a consequence of the vanishing
of those terms, and when $Q$ takes all the possible values, i.e.,
all the possible branches of the singular points on $\chi$ and of
the support of $R$, we get the linear equations globally imposed
by the condition (\ref{eq:ad2}). An easy reasoning reveals that the number of such adjoint
conditions is equal to
\begin{equation} \label{eq:nadj}
\frac{1}{2} \deg \ad + \deg R = \frac{1}{2} \sum_{P \in \chi} n_P
+ \deg R = \sum_{P \in \chi} \delta_P + \deg R.
\end{equation}

\begin{ejem}
Let $\chi$ be the projective plane curve over the finite field of two elements $\mathbb{F}_2$ given
by the equation $F(X,Y,Z)=X^3-Y^2Z$. The only singular point of
$\chi$ is $P_1=[0:0:1]$. Let be the point $P_2=[0:1:0]$ and the
effective divisor $R=0P_1+P_2$. The adjunction divisor of $\chi$
is $\ad=2 P_1$. A local equation of $\chi$ with $P_1=(0,0)$ is
$f(x,y)=x^3-y^2$. A parametrization of $\chi$ at $P_1$ is given
by
\begin{displaymath}
  \begin{array}{l}
  X_1(t_1)=t_1^2\\
  Y_1(t_1)=t_1^{3}
\end{array} 
\end{displaymath}
Take a form $H \in \mathcal{F}_{4-3=1}$, $H(X,Y,Z)=aX+bY+cZ$. First we
want to express the adjoint conditions in  terms of the coefficients
\[
\NNA H \ge \ad + R = 2P_1+P_2.
\]
To this end we consider first a local equation for $H$ at $P_1$, namely
\[
h(x,y)=H(x,y,1)= ax+by+c.
\]
Then $h(X_1(t_1),Y_1(t_1))=h(t^2,t^3)=at_1^2+bt_1^3+c$. So if we
wish to have
\[
\mathrm{ord}_{t_1}(h(X_1(t_1),Y_1(t_1)))=\mathrm{ord}_{t_1}(bt_1^3+at_1^2+c)
\ge 2
\]
(since $2$ is the coefficient for $P_1$ and
$(X_1(t_1),Y_1(t_1))$ is a parametrization at $P_1$), then this is
possible if and only if $c=0$. Thus $c=0$ is one of the required
\emph{linear} adjoint conditions. \vspace{0.25cm}

Now consider a local equation for $\chi$ at $P_2$. This is
$f^{\prime}(x,z)=F(x,1,z)=x^3-z$, and admits a parametrization
\begin{displaymath}
 \begin{array}{l}
  X_2(t_2)=t_2\\
  Z_2(t_2)=t_2^{3}
\end{array} 
\end{displaymath}
Consider the local equation for $H$ at $P_2$
$$
h^{\prime}(x,z)=H(x,1,z)=ax+b+cz.
$$
Hence the adjoint conditions imposed by $\NNA H \ge \ad + R=2 P_1
+ P_2$ at $P_2$ come from considering
$h^{\prime}(X_2(t_2),Z_2(t_2))=h^{\prime}(t_2,t_2^3)=at_2+b+ct_2^3$
and they impose
$$
\mathrm{ord}_{t_2}(h^{\prime}(X_2(t_2),Z_2(t_2)))=\mathrm{ord}_{t_2}(ct_2^3+at_2+b)
\ge 1.
$$
This inequality holds whenever $b=0$. Hence $b=0$ is another
\emph{linear} equation taking part in the set of adjoint
conditions contained in $\NNA H \ge \ad +R$. We have obtained
two adjoint conditions, as we had hoped by (\ref{eq:nadj}),
since $\frac{1}{2}\deg \ad + \deg R = \frac{1}{2} \cdot 2 + 1 =2$.
\end{ejem}

We conclude this section with two remarkable results. Let $\chi$ be an absolutely irreducible projective plane curve
defined over a perfect field $\ff$ and given by an equation
$F(X_0,X_1,X_2)=0 $, where $F \in \mathcal{F}_d$. One
application of the adjoint forms is the following result, due to
Max Noether (he stated it of course not in this way; our version
may be found in Hach\'e and Le Brigand \cite{hache}, Theorem 4.2, and
Le Brigand and Risler \cite{lebri}, \S 3.1):
\begin{teo}[Max Noether] Let $\chi, \chi^{\prime}$ be curves as
above given by homogeneous equations $F(X_0,X_1,X_2)=0$ and
$G(X_0,X_1,X_2)=0 $ respectively and such that $\chi^{\prime}$
does not contain $\chi$ as a component. Then, if we consider
another such a curve given by $H(X_0,X_1,X_2)=0$ with $\NNA H
\ge \ad + \NNA G$ (where $\ad$ is the adjunction divisor on
$\chi$), there exist forms $A,B$ with coefficients in $\ff$ such
that $H=AF+BG$.
\end{teo}

This theorem has great importance, and, for instance, allows us to
prove the Brill-Noether theorem, which 
gives a basis for the vector spaces $\mathcal{L}(D)$. Readers are referred to \cite{hache}, Theorem 4.4, for further details. A short remark about notation is needed.
For any non effective divisor $D$ we will write
$D=D_{+}-D_{-}$ with $D_{+}$ and $D_{-}$ effective divisors
of disjoint support. 

\begin{teo}[Brill-Noether] \label{teo:br}
Let $\chi$ be an adjoint curve as above with normalization $\wchi$. Let $\ad$ be its adjunction divisor and let $D$ be a divisor on
$\wchi$ rational over $\ff$. Moreover, consider a form $H_0 \in
\mathcal{F}_{\widetilde{n}}$ defined over $\ff$, not divisible by
$F$ and satisfying $\NNA H_0 \ge \ad + D_{+}$. Then
\[
\mathcal{L}(D)=\Bigg \{\frac{h}{h_0} \mid H \in
\mathcal{F}_{\widetilde{n}}, F \nmid H \ \mathrm{and} \ \NNA H + D
\ge \NNA H_0  \Bigg\} \cup \{0\},
\]
where $h,h_0 \in \ff (\chi)$ denote respectively the rational
functions $H,H_0$ restricted on $\chi$.
\end{teo}

\begin{rem} \label{rem:br}
Such a form $H_0 \in \mathcal{F}_{\widetilde{n}}$ exists whenever
$\widetilde{n} > \max \bigg \{d-1, \frac{d-3}{2}+\frac{\deg (\ad +
D_+)}{d} \bigg \}$
(see Hach\'e and Le Brigand \cite{hache} for details).
\end{rem}

\section{The Weierstra\ss~semigroup at several points}

Let $\chi$ be an absolutely irreducible
projective algebraic plane curve defined over a perfect field
$\ff$. Let $\underline{P}$ denote a set of $r$ different points
$P_1,\ldots,P_r$ on $\chi$. Furthermore, the perfect field $\ff$
must have cardinality greater or equal to $r$: $\sharp \ff \ge r$.
Let $\widetilde{\chi}$ be the normalization of $\chi$.
\medskip

Our purpose is to compute the dimensions of the so-called
\emph{Riemann-Roch quotients}:
\[
0 \le \dim_{\ff} \frac{\mathcal{L}(\mep)}{\mathcal{L}(\mepu)} \le
r
\]
by choosing functions in
$\mathcal{L}(\mep)=\mathcal{L}(m_1P_1+\ldots+m_rP_r)$ but not in
$\mathcal{L}(\mepu)=\mathcal{L}((m_1-1)P_1+ \ldots + (m_r-1)P_r)$,
that is, achieving at the $P_i$ poles of order $m_i$. We are going
to restrict to the case $m_i \in \mathbb{N}$, for all
$i=1,\ldots,r$. Such dimensions will be determined by the previous
calculus of the \emph{Riemann-Roch quotients with respect to
$P_i$}:
\[
0 \le \dim_{\ff}
\frac{\mathcal{L}(\mep)}{\mathcal{L}((\me-\ei)\underline{P})} \le
1,
\]
where $\ei$ denotes the vector in $\mathbb{N}^r$ with $1$ in the
$i$-th position and $0$ in the other ones. 
\medskip

Summarizing, this section deals with the following topics:
\begin{itemize}
\item[(a)] How to compute $\dim_{\ff}
\frac{\mathcal{L}(\mep)}{\mathcal{L}((\me-\ei)\underline{P})}$ and
an associated function belonging to this quotient vector space
when such a dimension is $1$.

\item[(b)] How to compute $\dim_{\ff}
\frac{\mathcal{L}(\mep)}{\mathcal{L}(\mepu)}$ (deducing bounds).

\item[(c)] How to compute the Weierstras\ss~semigroup at two
points.
\end{itemize}

All the statements and proofs of this section can be found in
\cite{cato}, \S 2.
\medskip

Consider a finite
set of nonsingular points $P_1, \ldots ,P_r$ on $\chi$ and a divisor
$m_1P_1+ \ldots +m_rP_r$ for $m_i \in \mathbb{N}$ $\forall \
i=1,\ldots r$. We will denote $\pe = \{P_1, \ldots , P_r \}$,
$\mep=m_1P_1+\ldots +m_rP_r$, $\me=(m_1,\ldots,m_r)$,
$\ei=(0,\ldots,0,1,0,\ldots,0)$ and $\uno=(1,\ldots,1)$.

\begin{defi}
For $\pe \in \chi$ we define 
\[
\Gamma_{\pe}:= \Big \{
-(\mathrm{ord}_{P_1}(f),\ldots,\mathrm{ord}_{P_r}(f)) \mid f \in
\ff(\chi)^{\ast} \ \mathrm{regular \ at} \ \chi \setminus \pe
\Big\}.
\]
\end{defi}

Obviously $\wss$ is a subsemigroup of $(\mathbb{N},+)$. 
Notice that, for $\mep = m_1 P_1+m_2 P_2$, the fact that $f \in
\mathcal{L}(\mep)$ is equivalent to the inequalities
\begin{displaymath}
(\star) \left\{ \begin{array}{l}
\mathrm{ord}_{P_1}(f) \ge -m_1\\
\mathrm{ord}_{P_2}(f) \ge -m_2.
\end{array} \right.
\end{displaymath}
This means: the set of possible orders ($\star$) which can be
taken by the function $f$ is represented by the shadowed area in the
figure (each axis represents one of the two branches):
\vspace{2cm}
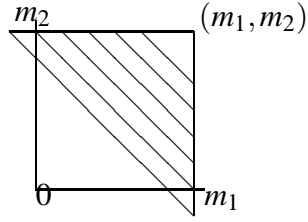
\begin{center}
 \setlength{\unitlength}{0.7cm}
\begin{picture}(-5,-3)
\put(-4,0){\line(1,0){3,2}}
\put(-4,0){\line(0,1){3,2}}
\put(-4.5,3){\line(1,0){3.5}}
\put(-1,-0.5){\line(0,1){3.5}}
\put(-1,0){\line(-1,1){3}}
\put(-1,.5){\line(-1,1){2.5}}
\put(-1,1){\line(-1,1){2}}
\put(-1,1.5){\line(-1,1){1.5}}
\put(-1,2){\line(-1,1){1}}
\put(-1,2.5){\line(-1,1){.5}}
\put(-1,-.5){\line(-1,1){3.5}}
\put(-4,-.3){0}
\put(-4.4,3.2){$m_2$}
\put(-.8,-.3){$m_1$}
\put(-.9,3.1){$(m_1,m_2)$}
\end{picture}
\end{center}
\vspace{0.6cm}

\begin{defi}
An element $\me \in \mathbb{N}^r$ is called a \emph{non-gap} of
$\pe$ if and only if $\me \in \wss$. Otherwise $\me$ is called a
\emph{gap}.
\end{defi}

A very important characterization for the non-gaps is given by the
following (see \cite{de}, p. 629):

\begin{lema} \label{lema:felix}
If $\me \in \mathbb{Z}^r$ then one has:
$$
\me \in \wss \ \ \ \mathrm{if \ and \ only \ if} \ \ \
\ell(\mep)=\ell((\me-\ei)\underline{P})+1 \ \forall \
i=1,\ldots,r.
$$
\end{lema}

For every $i=1,\ldots,r$ and $\me=(m_1,\ldots,m_r)
\in \mathbb{N}^r$, we set
\[
\nabla_i(\me):=\Big \{(n_1,\ldots,n_r) \in \wss \mid n_i=m_i \
\mathrm{and} \ n_j \le m_j \ \forall j \ne i \Big \}.
\]

Then the two conditions proven to be equivalent in Lemma \ref{lema:felix} are indeed also equivalent to $\nabla_i(\me)\neq 0$ for every $i \in \{1, \ldots , r\}$. 
\medskip

A gap $\me$ satisfying $\ell(\mep)=\ell((\me -\epsilon_i)\underline{P})$ for all $i \in \{1, \ldots , r\}$ (or, equivalently, such that $\nabla_i(\me)=\emptyset$ for all $i \in \{1, \ldots ,r\}$) is called \emph{pure}. It is easily seen: if $\me$ is a pure gap, then $m_i$ is a gap for $\Gamma_{P_i}$ for every $i\in \{1,\ldots ,r\}$. Furthermore, if $\underline{1} \in \Gamma_{\underline{P}}$, then there are no pure gaps. The converse does not hold, as Example  \ref{ejem:2wss} will show.
\medskip

A basic tool on Weierstra\ss~semigroups is the following

\begin{teo}[Weierstra\ss~gap~theorem]
Let $\wchi$ be a curve of genus $g\geq 1$. Let $P$ be a rational branch on $\wchi$. Then there are $g$ gaps $\gamma_1, \ldots , \gamma_g$ such that
\[
1=\gamma_1<\ldots <\gamma_g \leq 2g-1.
\]
\end{teo}


\begin{prop}
Let $\me=(m_1,\ldots,m_r) \in \mathbb{N}^r$. If $\me$ is a gap, then there exists a regular
differential form $\omega$ on $\wchi$ with $(\omega) \ge \me -
\ei$ and a zero at $P_i$ of order $m_i-1$ for \emph{some} $i \in
\{1,\ldots,r \}$.
\end{prop}

\dem~After Riemann-Roch theorem it is clear that
\begin{eqnarray*}
\ell (\mep) - i (\mep) &=& m_1+m_2+\ldots + m_r + 1 - g\\
\ell ((\me-\ei)\underline{P}) - i ((\me - \ei)\underline{P}) &=&
m_1+m_2+\ldots +m_r - 1 + 1 - g.
\end{eqnarray*}
By adding both equations we have
\[
[\underbrace{\ell(\mep)-\ell((\me-\ei)\underline{P})}_{=
\varphi(\me)}]-[\underbrace{i(\mep)-i((\me-\ei)\underline{P})}_{=
\psi(\me)}]=1,
\]
for every $i=1,\ldots,r$, where $0 \le \varphi(\me) \le 1$ and $-1 \le \psi (\me) \le 0$,
and therefore
\[
\varphi(\me)=1 \Leftrightarrow
\ell(\mep)-\ell((\me-\ei)\underline{P})=1 \Leftrightarrow
\ell(\mep)=\ell((\me-\ei)\underline{P})+1 \Leftrightarrow \me \in
\wss \Leftrightarrow \psi(\me)=0.
\]
Hence if $\me \notin \wss$ then $\dim_{\ff} \Big (
\frac{\Omega((\me-\ei)\underline{P})}{\Omega(\mep)} \Big)=1$ and
so there exists a regular differential form $\omega$ on $\wchi$
with $(\omega) \ge \me - \ei$ and
$\mathrm{ord}_{P_i}(\omega)=m_i-1$ for \emph{some} $i \in
\{1,\ldots,r \}$. \qed

\begin{prop}
Let $\chi$ be a plane curve of genus $g$, let $\pe$ be a set of $r$
closed points on $\chi$ and set
$\underline{m}=(m_1,\ldots,m_r) \in \mathbb{N}^r$. If $\me$ is a
gap, then $m_1+\ldots + m_r <2g$.
\end{prop}

\dem~Denote by $D_{2g,\pe}$ a divisor with degree $2g$ and support
$\pe$, and by $D_{2g-1,\pe}$ a divisor with degree $2g-1$ and
support $\pe$.
If $m_1+\ldots +m_r \ge 2g-1$ then $m_1+\ldots + m_r \ge 0$ as a consequence of 
Riemann-Roch, and for every
4$i=1,\ldots,r$
\[
\ell(D_{2g,\pe})=2g+1-g=g+1 \ne g=2g-1+1-g=\ell(D_{2g-1,\pe}),
\]
which implies that $\me$ is a non-gap, i.e., $\me \in \wss$. So,
if $\me \notin \wss$, then $m_1+\ldots +m_r < 2g$. \qed
\medskip

Notice that, for divisors of the form $\mep = m_1 P_1 +
m_2 P_2$, the plane $\mathbb{N} \times \mathbb{N}$ is divided in
three parts by the line $m_1+m_2=2g$ as in the figure, namely


\vspace{3.5cm}

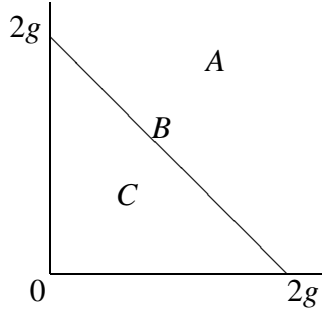
\begin{center}
 \setlength{\unitlength}{0.9cm}
\begin{picture}(-5,-3)
\put(-4,0){\line(1,0){4}}
\put(-4,0){\line(0,1){4}}
\put(-.5,0){\line(-1,1){3.5}}
\put(-4.3,-.4){0}
\put(-4.6,3.5){$2g$}
\put(-.5,-.4){$2g$}
\put(-1.7,3){$A$}
\put(-2.5,2){$B$}
\put(-3,1.01){$C$}
\end{picture}
\end{center}


\vspace{0.1cm}

$$
A:=\{(m_1,m_2) \mid m_1+m_2>2g, \ m_1>0, \ m_2>0 \}
$$
$$
B:=\{(m_1,m_2) \mid m_1+m_2=2g, \ m_1>0, \ m_2>0 \} \cup
\{(m_1,0), \ m_1>2g \} \cup \{(0,m_2), \ m_2>2g \}
$$
$$
C:=\{(m_1,m_2) \mid m_1+m_2<2g, \ m_1 \ge 0, \ m_2 \ge 0 \}.
$$

All the points lying on $A$ and $B$ correspond to values in $\wss$, but nothing can be a priory said
about the points on $C$.

\subsection{Dimension of the Riemann-Roch quotients with respect to $P_i$ and associated functions}

We start by computing the dimension of the Riemann-roch quotients associated to the points $P_i$.

\begin{prop} \label{prop:2}
Let $\me \in \mathbb{N}^r$ such that $\sum_{i=1}^r m_i <2g$. Then,
for $i \in \{1, \ldots, r\}$ we have:
\begin{itemize}
\item[a)] $\mathrm{dim}_{\ff} [\Omega((\me - \ei)\underline{P})
\setminus \Omega (\mep)]=1$ if and only if $\exists$ a homogeneous
polynomial $H_0$ of degree $d-3$ with $\NNA H_0 \ge \ad +
(\me-\ei)\pe$ such that $P_i$ is not in the support of the
effective divisor $\NNA H_0 - \ad - (\me-\ei)\pe$.

\item[b)] $\exists$ $\underline{m}^{\prime} \ge \me$ with
$\mathrm{dim}_{\ff} [\Omega((\underline{m}^{\prime} -
\ei)\underline{P}) \setminus \Omega (\underline{m}^{\prime}
\underline{P})]=1$ if and only if $\exists$ a homogeneous
polynomial $H_0$ of degree $d-3$ such that $\NNA H_0 \ge \ad +
(\me-\ei)\pe$.
\end{itemize}
\end{prop}

\dem~\begin{itemize}
\item[a)] If $\mathrm{dim}_{\ff} [\Omega((\underline{m} -
\ei)\underline{P}) \setminus \Omega (\underline{m}
\underline{P})]=1$, then this is equivalent to $\me \notin \wss$
and also to the existence of an index $i$ with
$\ell(\mep)=\ell((\me-\ei)\pe)$, or, in other words, to the existence of an index $i$
with $i((\me-\ei)\pe)=i(\mep)+1$; that is, there exists a
homogeneous polynomial $H_0$ of degree $d-3$ such that $\NNA H_0
\ge \ad + (\me-\ei)\pe$.

\item[b)] If there is $\me^{\prime} \ge \me$ with
$\mathrm{dim}_{\ff} [\Omega((\underline{m}^{\prime} -
\ei)\underline{P}) \setminus \Omega (\underline{m}^{\prime}
\underline{P})]=1$ then there exists an adjoint $H_0$ of degree
$d-3$ whose divisor is $\ge (\me^{\prime}-\ei)\pe$ outside $\ad$,
i.e., $\NNA H_0 - \ad \ge (\me^{\prime}-\ei)\pe \ge (\me-\ei)\pe$.
Conversely, if there is $H_0$ of degree $d-3$ with $\NNA H_0 \ge
\ad + (\me-\ei)\pe$ then there exists $\omega \ne 0$ differential
form such that $(\omega)=\NNA H_0-\ad \ge (\me-\ei)\pe$. Assume
that $\me^{\prime}-\ei$ are the orders of the zeros of $\omega$ at
$\pe$. Thus, $\me^{\prime}-\ei \ge \me-\ei$, what implies 
$\me^{\prime} \ge \me$ and $\omega \in \Omega((\me -
\ei)\underline{P}) \setminus \Omega (\mep)$.\hfill \qed
\end{itemize}

The following corollary yields a way to
relate the adjunction theory and the computation of the
Weierstra\ss~semigroup at several points:

\begin{cor} \label{cor:one}
Let $\me \in \mathbb{N}^r$ with $\sum_{i=1}^{r}m_i <2g$. For a
given form $H$ of degree $d-3$ and $i \in \{1,\ldots,r\}$ there
exists a condition imposed by the inequality $\NNA H \ge \ad +
\mep$ at $P_i$ which is independent of the conditions imposed by
$\NNA H \ge \ad +(\me - \ei)\underline{P}$ at $P_i$ if and only if
\[
\dim_{\ff} \frac{\Omega ((\me - \ei)\underline{P})}{\Omega
(\mep)}=1.
\]
\end{cor}

The second step is the computation of the rational functions associated to the nongaps of the Weierstra\ss~semigroup.
Note that, if $\dim_{\ff}\frac{\Omega ((\me-\ei)\underline{P})}{\Omega
(\mep)}=0$, then $\dim_{\ff}\frac{\mathcal{L}(\mep)}{\mathcal{L}
((\me-\ei)\underline{P})}=1$ and so there is a rational function $f_{i,\me} \in
\frac{\mathcal{L}(\mep)}{\mathcal{L} ((\me-\ei)\underline{P})}$
with a pole of order $m_i$ at $P_i$. In order to compute such a
function, we base on Brill-Noether Theorem \ref{teo:br}:

\begin{algo} \label{algo:no}
Preserving notations as above, we obtain a function $f_{i,\me} \in
\frac{\mathcal{L}(\mep)}{\mathcal{L} ((\me-\ei)\underline{P})}$
with a pole of order $m_i$ at $P_i$ by following these steps:

\begin{itemize}
\item[-] Compute a homogeneous polynomial $H_0$ not divisible by
$F$ of large enough degree $n$ in the sense of Remark \ref{rem:br}
satisfying $\NNA H_0 \ge \ad + \mep$.

\item[-] Calculate $R_{\me}$, which is the effective divisor such
that $\NNA H_0 = \ad + \mep + R_{\me}$. Obviously
$R_{\me-\ei}=R_{\me}+P_i$.

\item[-] Find a form $H_{\me}$ of degree $n$ not divisible by $F$
such that $\NNA H_{\me} \ge R_{\me}$ but not satisfying $\NNA
H_{\me} \ge R_{\me - \ei}=R_{\me}+P_i$.

\item[-] Output: $f_{i,\me}=\frac{h_{\me}}{h_0}$, where
$h_{\me},h_0$ are the restricted forms on $\chi$ for $H_{\me}$ and
$H_0$ respectively.
\end{itemize}
\end{algo}
\vspace{0.2cm}

\begin{ejem}
Let $\chi$ be the curve given by the equation
$F(X,Y,Z)=X^3Z+X^4+Y^3Z+YZ^3$ and consider the points
$P_1=[0:1:1]$ and $P_2=[0:1:0]$ and $\me=(1,2)$. We want to
compute $\dim_{\ff}\frac{\mathcal{L}(\mep)}{\mathcal{L}((\me -
\eu)\underline{P})}$ and
$\dim_{\ff}\frac{\mathcal{L}(\mep)}{\mathcal{L}(\me-\ed)\underline{P}}$.
\vspace{0.25cm}

A local parametrization of $F$ at $P_1$ is given by
\begin{displaymath}
   \begin{array}{l}
  X_1(t_1)=t_1\\
  Y_1(t_1)=t_1^{3}+t_1^{4}+t_1^9+t_1^{10}+t_1^{11}+t_1^{12}+\ldots
\end{array} 
\end{displaymath}

with local equation $f_1(x,y)=y^2+y^3+x^3+x^4$. Analogously at $P_2$
\begin{displaymath}
  \begin{array}{l}
  X_2(t_2)=t_2\\
  Z_2(t_2)=t_2^{4}+t_2^7+t_2^{10}+t_2^{12}+t_2^{13}+t_2^{16}+\ldots
\end{array} 
\end{displaymath}
with local equation $f_2(x,z)=z+z^3+x^3z+x^4$. \vspace{0.25cm}

First, we calculate the adjunction divisor: this is $\ad=2P_1$.
\vspace{0.25cm}

Search a form $H$ of degree $d-3=4-3=1$, that is, a linear form
$H(X,Y,Z)=aX+bY+cZ$. At $P_1$ $H$ admits the equation
$h_1(x,y)=H(X,Y-1,1)=ax+by+b+c$. At $P_2$ $H$ admits the equation
$h_2(x,z)=H(X,1,Z)=ax+b+cz$. Then:

\begin{displaymath}
 \begin{array}{l}
h_1(X_1(t_1),Y_1(t_1))=at_1+b(t_1^3+t_1^4+t_1^9+\ldots)+b+c=(b+c)+at_1+bt_1^3+bt_1^4+bt_1^9+\ldots.\\
h_2(X_2(t_2),Z_2(t_2))=b+at_2+ct_2^4 +ct_2^7+ct_2^{10}+\ldots.
\end{array}
\end{displaymath}
\vspace{0.2cm}

In order to compute
$\dim_{\ff}\frac{\mathcal{L}(\mep)}{\mathcal{L}((\me -
\eu)\underline{P})}$ we impose the systems of equations with the
adjunction conditions at $P_1$: \vspace{0.25cm}

\begin{displaymath}
  \left\{ \begin{array}{l}
  \NNA H \ge \ad + (m_1-1)P_1=2P_1\\
  \NNA H \ge \ad + m_1 P_1 = 3P_1,
\end{array} \right.
\end{displaymath}
or, in other words

\begin{displaymath}
  \left\{ \begin{array}{l}
  \mathrm{ord}_{t_1}(h_1(X_1(t_1),Y_1(t_1))) \ge 2 \Longrightarrow b+c=a=0\\
  \mathrm{ord}_{t_1}(h_1(X_1(t_1),Y_1(t_1))) \ge 3 \Longrightarrow b+c=a=0
\end{array} \right.
\end{displaymath}
\vspace{0.2cm}

So the second system does not add any independent condition to the
first one; this means, by Corollary \ref{cor:one}, that
$\dim_{\ff}\frac{\mathcal{L}(\mep)}{\mathcal{L}((\me -
\eu)\underline{P})}=1$. \vspace{0.25cm}

In order to compute
$\dim_{\ff}\frac{\mathcal{L}(\mep)}{\mathcal{L}((\me -
\ed)\underline{P})}$ the systems of equations with the adjunction
conditions at $P_2$ are
\begin{displaymath}
  \left\{ \begin{array}{l}
  \NNA H \ge (m_2-1)P_2=P_2\\
  \NNA H \ge m_2 P_2=2P_2,
\end{array} \right.
\end{displaymath}
that is,
\begin{displaymath}
  \left\{ \begin{array}{l}
  \mathrm{ord}_{t_2}(h_2(X_2(t_2),Z_2(t_2))) \ge 1 \Rightarrow
  b=0\\
  \mathrm{ord}_{t_2}(h_2(X_2(t_2),Z_2(t_2))) \ge 2 \Rightarrow a=0=b.
\end{array} \right.
\end{displaymath}
\vspace{0.2cm}

 Notice that, in this case, the adjunction divisor does not appear in the inequalities since $P_2$ does not belong to its support.
 The second system adds one independent condition to the first one
  and this means that
$\dim_{\ff}\frac{\mathcal{L}(\mep)}{\mathcal{L}((\me -
\ed)\underline{P})}=0$ again by Corollary \ref{cor:one}. \qed
\end{ejem}
\vspace{0.2cm}

\begin{ejem}
Consider the previous example but with $\me=(4,6)$. As
$m_1+m_2=4+6=10 >2g$, we know without calculations $\me \in
\wss$, i.e., that
$\dim_{\ff}\frac{\mathcal{L}(\mep)}{\mathcal{L}((\me -
\ei)\underline{P})}=1$ for $i=1,2$. So we will look for the
corresponding functions $f_{i,\me}$ with poles at $P_i$ of order
$m_i$ for $i=1,2$. \vspace{0.25cm}

First of all, we search $\widetilde{n}>\max \Big
\{3,\frac{2}{4}+\frac{12}{4} \Big \}=\max \Big\{3,\frac{14}{4}
\Big\}$. Let us take $\widetilde{n}=5$. \vspace{0.25cm}

Then we look for a form $H_0$ of degree $\widetilde{n}=5$ such
that $\NNA H_0 \ge \ad + \mep$. In this case $\NNA H_0 \ge
4P_1+6P_2+2P_3$, since $\ad=2P_3$, where $P_3=[0:0:1]$. After some computations we find $H_0=X^4Z$. \vspace{0.25cm}

In order to compute $\NNA H_0$, we have to calculate $\NNA (X)$,
$\NNA (Y)$ and $\NNA (Z)$. Intersection points between $\{F=0\}$
and $\{X=0\}$ are $P_1=[0:1:1]$, $P_2=[0:1:0]$ and $P_3=[0:0:1]$
with multiplicities $1$, $1$ and $2$ respectively. So $\NNA
(X)=P_1+P_2+2P_3$. Intersection points between $\{F=0\}$ and
$\{Y=0\}$ are $P_3=[0:0:1]$ and $P_4=[1:0:1]$ such that
$\NNA (Y)=3P_3+P_4$. And the only point lying in the intersection
between $\{F=0 \}$ and $\{Z=0 \}$ is $P_2=[0:1:0]$ with
multiplicity $4$, therefore $\NNA (Z)=4P_2$. \vspace{0.25cm}

Thus $\NNA H_0 = 4 \NNA (X)+\NNA (Z)=4P_1+8P_2+8P_3$. The residue
divisor $R_{\me}=\NNA H_0 - \ad -\mep = 2P_2+6P_3$. Following the
algorithm described above, we have to find a form $H_{\eu}$ such
that $\NNA H_{\eu} \ge R_{\me}$ but $\NNA H_{\eu} \ngeq R_{\me} +
P_1$. For instance we take $H_{\eu}=Y^2Z^3$, since

\begin{displaymath}
 \begin{array}{l}
\NNA H_{\eu}=12P_2+6P_3+2P_4 \ge 2P_2 +6P_3\\
\NNA H_{\eu} \ngeq P_1+2P_2+6P_3.
\end{array}
\end{displaymath}
\vspace{0.2cm}

So $f_{1,\me}=\frac{Y^2Z^3}{X^4Z}=\frac{Y^2Z^2}{X^4} \in
\frac{\mathcal{L}(\mep)}{\mathcal{L}((\me-\eu)\underline{P})}$.
\vspace{0.25cm}

Now we have to find a form $H_{\ed}$ such that $\NNA H_{\ed} \ge
R_{\me}$ but $\NNA H_{\ed} \ngeq R_{\me} + P_2$. For instance we
take $H_{\ed}=X^2Y^3$, since

\begin{displaymath}
 \begin{array}{l}
\NNA H_{\ed}=2P_1+2P_2+13P_3+3P_4 \ge 2P_2 +6P_3\\
\NNA H_{\ed} \ngeq 3P_2+6P_3.
\end{array}
\end{displaymath}
\vspace{0.2cm}

Thus $f_{2,\me}=\frac{X^2Y^3}{X^4Z}=\frac{Y^3}{X^2Z} \in
\frac{\mathcal{L}(\mep)}{\mathcal{L}((\me-\ed)\underline{P})}$.\qed
\end{ejem}

\begin{algo} \label{algo:si}
There is an alternative way of calculating these functions
$f_{i,\me}$, computationally more effective:

\begin{enumerate}
\item Take a basis of $\mathcal{L}(\mep)$, say $\{h_1, \ldots
,h_s\}$.
\item Calculate the pole orders at $P_i$,
$\{-\mathrm{ord}_{P_i}(h_1),\ldots, -\mathrm{ord}_{P_i}(h_s) \}$.
\item Order these pole orders increasing, in such a way that
$-\mathrm{ord}_{P_i}(h_s)=m_i$. We can assume this, as otherwise,
if $-\mathrm{ord}_{P_i}(h_s)=k_i > m_i$ we can replace $m_i$ by
$k_i$, since $\mathcal{L}(m_1 P_1 + \ldots + m_i P_i + \ldots +m_r
P_r)=\mathcal{L}(m_1 P_1 + \ldots + k_i P_i + \ldots + m_r P_r)$.
\item The function $h_s$ has pole order $m_i$ at $P_i$, but other
functions could also have the same property. So, for any $h_j$
satisfying $-\mathrm{ord}_{P_i}(h_j)=m_i$, there exists $\lambda_j
\ne 0$ in $\ff$ such that $h_j=\lambda_j h_s$, that is,
$-\mathrm{ord}_{P_i}(h_j-\lambda_j h_s)<m_i$. So we change $h_j$
by $g_j:=h_j-\lambda_j h_s$, and $g_k:=h_k$ for $k \ne j$.
\item Now we have a set of functions $\{g_1, \ldots, g_s \}$ where
$g_s=f_{i,\me}$, and $\{g_1,\ldots,g_{s-1}\}$ is a basis
for the vector space $\mathcal{L}((\me-\ei)\underline{P})$.
\end{enumerate}
\end{algo}

\begin{ejem}
We present a worked example in \textsc{Singular} for computing functions as
above. First we import the library \verb"brnoeth.lib" and another
one \verb"several.lib" in which we have programmed the procedure
\verb"ordRF" that computes the pole orders of a rational
function:{\small
\begin{verbatim}
> LIB "brnoeth.lib";
> LIB "several.lib";
> int plevel=printlevel;
> printlevel=-1;
\end{verbatim}}
We define the ring and the curve:{\small
\begin{verbatim}
> ring s=2,(x,y),lp;
> list C=Adj_div(x3y+y3+x);
 ==>The genus of the curve is 3
\end{verbatim}}
 The list of computed places is{\small
\begin{verbatim}
> C=NSplaces(1,C);
> C[3];
 -->[1]:
 -->   1,1
 -->[2]:
 -->   1,2
 -->[3]:
 -->   1,3
\end{verbatim}}
The base point of the first place of degree $1$
is, in homogeneous coordinates:{\small
\begin{verbatim}
> def SS=C[5][1][1];
> setring SS;
> POINTS[1];
 -->[1]:
 -->   0
 -->[2]:
 -->   1
 -->[3]:
 -->   0
> setring s;
\end{verbatim}}
We define the divisor \verb"G=4C[3][1]+4C[3][3]":{\small
\begin{verbatim}
> intvec G=4,0,4;
> def R=C[1][2];
\end{verbatim}}
A basis \verb"LG" for $\mathcal{L}(\mep)$ is
supplied by the Brill-Noether algorithm:{\small
\begin{verbatim}
> setring R;
> list LG=BrillNoether(G,C);
 -->Vector basis successfully computed
> int lG=size(LG);
\end{verbatim}}
The pole orders for the rational functions in \verb"LG"
are{\small
\begin{verbatim}
> int j;
> intvec h;
> for (j=1;j<=lG;j=j+1){
. h[j]=ordRF(LG[j],SS,1)[1]; . }
> h;
 -->0,-1,2,-2,-3,-4
\end{verbatim}}
 And the desired rational function is{\small
\begin{verbatim}
> LG[lG];
 -->_[1]=xyz2+y4
 -->_[2]=x4
> printlevel=plevel;
\end{verbatim}}
\end{ejem}

\subsection{Dimension of the Riemann-Roch quotients}

Computing the dimension of
$\frac{\mathcal{L}(\mep)}{\mathcal{L}(\mepu)}$ is an easy task
by Corollary \ref{cor:one}:
\begin{prop}
$$
\frac{\mathcal{L}(\mep)}{\mathcal{L}(\mepu)}=\frac{\mathcal{L}(\mep)}{\mathcal{L}((\me-\eu)\underline{P})}
\oplus
\frac{\mathcal{L}((\me-\eu)\underline{P})}{\mathcal{L}((\me-\eu-\ed)\underline{P})}
\oplus \ldots \oplus \frac{\mathcal{L}((\me-\eu-\ed-\ldots
-\er)\underline{P})}{\mathcal{L}(\mepu)}.
$$
\end{prop}

\dem~It is just to define the map $(f_1,\ldots,f_r) \mapsto f_1 + \ldots + f_r$.
\qed

\begin{nota}
Notice that the map $(f_1,\ldots,f_r) \mapsto f_1 \cdot
\ldots \cdot f_r$ cannot work, because products do not preserve poles.
It is also important the fact that the $P_i$ are
different, otherwise the statement does not hold: take for example
$f(z)=\frac{1}{z}$ and $h(z)=z-\frac{1}{z}$ in $\mathbb{C}$.
The sum $f(z)+h(z)=\frac{1}{z}+z-\frac{1}{z}=z$ has no poles,
however $f(z)$ and $h(z)$ have both a simple pole at $0$.
\vspace{0.25cm}
\end{nota}

\subsection{Computing the Weierstra\ss~semigroup at several points}

Preserving notations, let $\wss$ be the Weierstra\ss~semigroup at
the points $P_1,P_2, \ldots, P_r$ and $\Gamma_{P_i}$ the
Weierstra\ss~semigroups corresponding to the points $P_i$ for
$i=1, \ldots, r$. Write $\mathbb{N}^{\ast}:=\mathbb{N} \setminus
\{ 0 \}$ and $\me_i:=\me - m_i \ei$. 

\begin{prop}\label{cor:211}
Let $\me \in \mathbb{N}^r$, $i \in \{1,\ldots,r \}$ and $\me_i \in
\mathbb{N}^r \setminus \wss$. Let
$$
m:=\mathrm{min}\Big \{ n \in \mathbb{N}^{\ast} \mid \me_i+n \ei
\in \wss \Big \}.
$$
Then any vector $\underline{n}=(n_1,\ldots,n_r) \in \mathbb{N}^r$
belongs to $\mathbb{N}^r \setminus \wss$ whenever $n_i=m$, and
$n_j=m_j=0$ or $n_j < m_j$ for $j \ne i$. In particular, $m$ is a
gap at $P_i$.
\end{prop}

Define the usual partial order $\preceq$ over
$\mathbb{N}^r$, that is, for $\me,\underline{n} \in \mathbb{N}^r$:
$$
(m_1,\ldots,m_r) \preceq (n_1,\ldots,n_r) \ \ \Longleftrightarrow
\ \ m_i \le n_i \ \mathrm{for \ all} \ i=1,\ldots,r.
$$

\begin{prop} \label{cor:main}
Let $i \in \{1, \ldots, r \}$ and $\me =(m_1,\ldots,m_r)$ be a
minimal element of the set
$$
\Big \{(n_1,\ldots,n_r) \in \wss \mid n_i=m_i \Big \}
$$
with respect to the partial order $\preceq$. Assume that $n_i >0$
and the existence of one $j \in \{1,\ldots,r \}$, $j \ne i$ with
$m_j >0$. Then:
\begin{itemize}
    \item[a)] $\me_i \in \mathbb{N}^r \setminus \wss$.
    \item[b)] $m_i=\mathrm{min}\{n \in \mathbb{N}^{\ast} \mid \me_i+n \ei \in \wss \}$;
    in particular, $m_i$ is a gap at $P_i$.
\end{itemize}
\end{prop}

Propositions \ref{cor:211} and \ref{cor:main} determine a surjective map 
\begin{displaymath}
\begin{array}{lccc}
\varphi_i: &\Big \{\me_i \in \mathbb{N}^r \mid \me_i \in \mathbb{N}^r \setminus \wss  \Big \}& \longrightarrow &\mathbb{N} \setminus \Gamma_{P_i}\\
& \me_i & \mapsto &  \min \Big \{m \in \mathbb{N}^{\ast} \mid
\me_i + m \ei \in \wss \Big \}.
\end{array}
\end{displaymath}

For $r=2$ this is in fact a bijection between the set of gaps at
$P_1$ and the set of gaps at $P_2$:
$$
m_1 \in \mathbb{N} \setminus \Gamma_{P_1} \Leftrightarrow (m_1,0)
\in \mathbb{N}^2 \setminus \wss \mapsto \beta_{m_1}:=\varphi_2
((m_1,0)) \in \mathbb{N}^ \setminus \Gamma_{P_2}.
$$
Furthermore, $m_1 = \min \Big \{n \in \mathbb{N}^{\ast} \mid
(n,\beta_{m_1}) \in \wss \Big \}$. More details can be found
in Homma and Kim \cite{hk} and Kim \cite{kim}. \vspace{0.25cm}

We summarize some remarkable facts for
the case of two points ($r=2$), which will be useful from
the computational point of view:
\begin{itemize}
    \item[(i)] All the gaps at $P_1$ and at $P_2$ are also gaps at
    $P_1,P_2$.
    \item[(ii)] By the Corollary \ref{cor:main}, for any gap $m_1$ at $P_1$, one has that
    $(m_1,\beta_{m_1})$ are gaps at $P_1,P_2$ for
    $\beta_{m_1}=0,1,\ldots,l_{m_1}$, until certain $0 \le
    l_{m_1}\le 2g-1$, with $g$ the genus of the curve and where
    $l_{m_1}$ satisfy that $l_{m_1}+1$ is a gap at $P_2$. The
    point $(m_1,l_{m_1}+1)$ is an element of $\wss$, which we will
    call the \emph{minimal (non-gap) at $m_1$}. We will refer to
    the set of the minimal non-gaps at every gap at $P_1$ (they will be $g$, since the number of gaps at $P_1$ is precisely $g$) as the
    set of \emph{minimal non-gaps at $P_1$}.
    \item[(iii)] The gaps obtained of that form, this is, the set 
     $$
    \Big \{(m_1,\beta_{m_1}) \in \mathbb{N}^2 \setminus \wss \mid m_1 \in \mathbb{N} \setminus \Gamma_{P_1} \ \mathrm{and} \ \beta_{m_1}=0,1,\ldots,l_{m_1} \ \mathrm{with} \ l_{m_1}+1 \in \mathbb{N} \setminus \Gamma_{P_2}  \Big \}
    $$
    will be called the \emph{set of gaps \emph{respect to} $P_1$}.
    \item[(iv)] Similarly, for any gap $m_2$ at $P_2$, one has that
    $(\alpha_{m_2},m_2)$ are gaps at $P_1,P_2$ for
    $\alpha_{m_2}=0,1,\ldots,l_{m_2}$, until some $0 \le
    l_{m_2}\le 2g-1$, with $g$ being the genus of the curve and where
    $l_{m_2}$ satisfy that $l_{m_2}+1$ is a gap at $P_1$. The
    point $(l_{m_2}+1,m_2)$ is an element of $\wss$, which we will
    call the \emph{minimal (non-gap) at $m_2$}. The set of the
    minimal non-gaps for every gap at $P_2$ will be called the set
    of \emph{minimal non-gaps at $P_2$}. The cardinality of such a set
    is $g$, since $g$ is the number of gaps at $P_2$.
    \item[(v)] The set of gaps
    $$
    \Big \{(\alpha_{m_2},m_2) \in \mathbb{N}^2 \setminus \wss \mid m_2 \in \mathbb{N} \setminus \Gamma_{P_2} \ \mathrm{and} \ \alpha_{m_2}=0,1,\ldots,l_{m_2} \ \mathrm{with} \ l_{m_2}+1 \in \mathbb{N} \setminus \Gamma_{P_1}  \Big \}
    $$
    is called the \emph{set of gaps \emph{respect to} $P_2$}.
    \item[(vi)] The intersection between the set of gaps respect to
    $P_1$ and respect to $P_2$ is not necessarily empty. In fact,
    the gaps in the intersection are just the \emph{pure gaps} at
    $P_1,P_2$.
\end{itemize}

The minimal
non-gaps at $P_1$ and $P_2$ provide enough information in order to deduce the
Weiestra\ss~semigroup at $P_1,P_2$. Recall that we have already
described algorithms to compute the dimension (and associated
functions, when is possible) of the Riemann-Roch quotients
$\frac{\mathcal{L}(\mep)}{\mathcal{L}((\me-\ei)\underline{P})}$
for given $\me$, $i \in \{1,2 \}$ and two rational points $P_1$,
$P_2$ on an absolutely irreducible projective algebraic plane
curve $\chi$ (see Algorithm \ref{algo:no} and Algorithm
\ref{algo:si}). An algorithm computing the set of minimal non-gaps
at $P_i$, for $i=1,2$ is the following:

\begin{algo}
Write $\dim(\me,P,C,i)$ for the
procedure calculating the dimension of the quotient vector space
$\frac{\mathcal{L}(\mep)}{\mathcal{L}((\me-\ei)\underline{P})}$:

{\small
\emph{\textsf{INPUT}}: points $P_1,P_2$, an integer $i \in \{
1,2\}$ and a curve $\chi$.\vspace{0.25cm}

\emph{\textsf{OUTPUT}}: the set of minimal non-gaps at
$P_i$.\vspace{0.25cm}

\begin{itemize}
\item let $L$ be a empty list and $g$ be the genus of $\chi$;

\item let $W_1$ and $W_2$  be the lists of gaps of $\chi$ at $P_1$
and $P_2$, respectively;

\item \verb"FOR" $k=1,\ldots,g$; $k=k+1$;
\begin{itemize}
\item \verb"IF" $i=1$ \verb"THEN"
\begin{itemize}
\item j=size of $W_2$;

\item \verb"WHILE" $\Big (\dim((W_1[k],W_2[j]),P,\chi,i)=1 \
\mathrm{AND} \  \dim((W_1[k],W_2[j]-1),P,\chi,i)=1) \ \mathrm{OR}
\ j=0 \Big )$ \verb"DO"
\begin{itemize}
\item $j=j-1$;
\end{itemize}
\item $L=L \cup \{(W_1[k],W_2[j]) \}$;

\item $W_2=W_2 \setminus \{j \}$;
\end{itemize}
\item \verb"ELSE"
\begin{itemize}
\item j=size of $W_1$;

\item \verb"WHILE" $\Big (\dim((W_1[j],W_2[k]),P,\chi,i)=1 \
\mathrm{AND} \  \dim((W_1[j]-1,W_2[k]),P,\chi,i)=1) \ \mathrm{OR}
\ j=0 \Big )$ \verb"DO"
\begin{itemize}
\item $j=j-1$;
\end{itemize}
\item $L=L \cup \{(W_1[j],W_2[k]) \}$;

\item $W_1=W_1 \setminus \{j \}$;
\end{itemize}
\end{itemize}
\item \verb"RETURN"($L$);
\end{itemize}
}
\end{algo}

\begin{ejem} \label{ejem:2wss}

Let $\chi$ be the curve over $\ff_2$ given by the equation
$F(X,Y,Z)=X^3Z+X^4+Y^3Z+YZ^3$. Consider the points
$P_1=[0:1:1]$ and $P_2=[0:1:0]$ on $\chi$. Then
\[
\mathbb{N}^2 \setminus \Gamma_{\{P_1,P_2\}}= \Big
\{(0,1),(0,2),(1,0),(1,2),(2,0),(2,1) \Big \},
\]
as shown in the figure (the black points are the elements of
$\wss$, the other ones are the gaps at $P_1,P_2$):

\vspace{4cm}

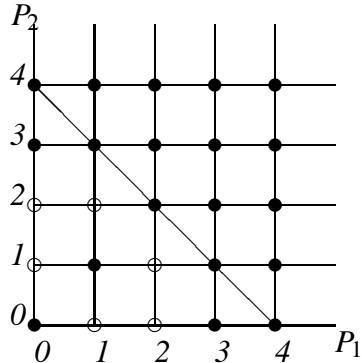
\begin{center}
\setlength{\unitlength}{1.6cm}
\begin{picture}(-5,-3)

\put(-4,0){\line(1,0){2.5}}

\put(-4,.5){\line(1,0){2.5}}

\put(-4,1){\line(1,0){2.5}}

\put(-4,1.5){\line(1,0){2.5}}

\put(-4,2){\line(1,0){2.5}}

\put(-4,0){\line(0,1){2.5}}

\put(-3.5,0){\line(0,1){2.5}}

\put(-3,0){\line(0,1){2.5}}

\put(-2.5,0){\line(0,1){2.5}}

\put(-2,0){\line(0,1){2.5}}

\put(-2,0){\line(-1,1){2}}

\put(-3.5,1.5){\circle*{.1}}

\put(-3.5,1){\circle{.1}}

\put(-3,0){\circle{.1}}

\put(-2,0){\circle*{.1}}

\put(-3,0.5){\circle{.1}}

\put(-2.5,2){\circle*{.1}}

\put(-2,1.5){\circle*{.1}}

\put(-2,1){\circle*{.1}}

\put(-3,1.5){\circle*{.1}}

\put(-4,1.5){\circle*{.1}}

\put(-4,2){\circle*{.1}}

\put(-4,1){\circle{.1}}

\put(-4,0.5){\circle{.1}}

\put(-4,0){\circle*{.1}}

\put(-3.5,2){\circle*{.1}}

\put(-3.5,0.5){\circle*{.1}}

\put(-3.5,0){\circle{.1}}

\put(-3,1){\circle*{.1}}

\put(-2.5,1){\circle*{.1}}

\put(-2.5,0.5){\circle*{.1}}

\put(-2.5,0){\circle*{.1}}

\put(-2,0.5){\circle*{.1}}

\put(-2,2){\circle*{.1}}

\put(-3,2){\circle*{.1}}

\put(-2.5,1.5){\circle*{.1}}

\put(-2,1){\circle*{.1}}

\put(-4.2,0){0}

\put(-4.2,.5){1}

\put(-4.2,1){2}

\put(-4.2,1.5){3}

\put(-4.2,2){4}

\put(-4.19,2.5){$P_2$}

\put(-4,-.3){0}

\put(-3.5,-.3){1}

\put(-3,-.3){2}

\put(-2.5,-.3){3}

\put(-2,-.3){4}

\put(-1.5,-.22){$P_1$}

\end{picture}
\end{center}
\vspace{1cm}

As an illustration of the Corollary \ref{cor:main}, for instance
let $i=1$, $\me=(m_1,m_2)=(2,2) \in \wss$ and the set $\Big \{
(n_1,n_2) \in \wss \mid n_1=m_1 \Big \}= \Big \{(2,n) \
\mathrm{for} \ n \ge 2 \Big \}$. A minimal element for this set is
$(2,2)$, and
\[
\me_i=\me - m_1 \eu = (2,2)-2(1,0)=(0,2)
\]
is a gap at $P_1,P_2$. We compute
\[
\min \Big \{n \in \mathbb{N}^{\ast} \mid (n,2) \in \wss  \Big
\}=2=m_1,
\]
and $m_1=2$ is actually a gap at $P_1$. \vspace{0.25cm}

In this example we can also see the bijection between the gaps at
$P_1$ and the gaps at $P_2$. Preserving notations as above, take
now $n_1=1$ as a gap at $P_1$. Then $(1,0)$ is a gap at $P_1,P_2$
and
\[
\varphi_2((1,0))=\min \Big \{n \in \mathbb{N}^{\ast} \mid
(1,0)+(0,n) \in \wss \Big \}=\min \Big \{n \in \mathbb{N}^{\ast}
\mid (1,n) \in \wss \Big \}=1,
\]
with $1$ being a gap at $P_2$. Moreover, $n_1=1=\min \Big \{n \in
\mathbb{N}^{\ast} \mid (n,1=\varphi_2((1,0))) \in \wss  \Big \}$.
\vspace{0.25cm}

Now take $p_1=2$ as the other gap at $P_1$. Then
$\varphi_2((2,0))=2$, which is a gap at $P_2$. Indeed $p_1=2=\min
\Big \{n \in \mathbb{N}^{\ast} \mid (n, \varphi_2((2,0)) \in \wss
\Big \}$. The same happens to the gaps at $P_2$.
\end{ejem}

\section{\textsf{Computational aspects using \textsc{Singular}}}

We are interested in explaining the most
important procedures implemented in \textsc{Singular} and to give
examples to show how to work with them. \vspace{0.25cm}

More precisely, in subsection \ref{section:hints}) we give some
hints of use of the library \verb"brnoeth.lib", since our
procedures are based on most of the algorithms contained in it.
Then, in Subsection \ref{section:seve} we present the procedures which
pretend generalize the computation of the Weierstra\ss~semigroup
to the case of several points, i.e., a set of procedures which
try to:
\begin{itemize}
\item[-] compute $\dim_{\ff}
\frac{\mathcal{L}(\mep)}{\mathcal{L}((\me-\ei)\underline{P})}$ and
a function $f_{\me,i} \in \mathcal{L}(\mep) \setminus
\mathcal{L}((\me-\ei)\underline{P})$ if possible.

\item[-] compute the set of minimal non-gaps at a point $P_i$, for
$i \in \{ 1,2 \}$.
\end{itemize}

\subsection{Hints of usage of \textsf{brnoeth.lib}}\label{section:hints}

The purpose of the library \verb"brnoeth.lib" of \textsc{Singular}
is the implementation of the Brill-Noether algorithm for solving
the Riemann-Roch problem and some applications in Algebraic
Geometry codes, involving the computation of
Weierstra\ss~semigroups for one point. \vspace{0.25cm}

A first warning: \verb"brnoeth.lib" accepts only prime base fields and
absolutely irreducible planes curves, although this is not
checked. \vspace{0.25cm}

Curves are usually defined by means of polynomials in two
variables, that is, by its local equation. It is possible to
compute most of the concepts concerning to the curve with the
procedure \verb"Adj_div". We defined the procedure (previously we
must have defined the ring, the polynomial $f$ and have charged
the library \verb"brnoeth.lib"):
\begin{verbatim}
> list C=Adj_div(f);
\end{verbatim}

The output consist of a list of lists as follows:
\begin{itemize}
    \item The first list contains the affine and the local ring.
    \item The second list has the degree and the genus of the
    curve.
    \item Each entry of the third list corresponds to one closed
    place,that is, a place and all its conjugates, which is
    represented by two integer, the first one the degree of the
    point and the second one indexing the conjugate point.
    \item The fourth one has the conductor of the curve.
    \item The fifth list consists of a list of lists, the first
    one, namely \verb"C[5][d][1]" being a (local) ring over an extension
    of degree $d$ and the second one (\verb"C[5][d][2]") containing the
    degrees of base points of places of degree $d$.
\end{itemize}

Furthermore, inside the local ring \verb"C[5][d][1]" we can find
the following lists:
\begin{itemize}
    \item \verb"list POINTS": base points of the places of degree $d$.
    \item \verb"list LOC_EQS": local equations of the curve at the base
    points.
    \item \verb"list BRANCHES": Hamburger-Noether expressions of the
    places.
    \item \verb"list PARAMETRIZATIONS": local parametrizations of the
    places.
\end{itemize}

Now we explain how the different kinds of common objects must be
treated in \verb"Singular". 
\medskip

\textbf{Affine points} $P$ are represented by a standard basis of
a prime ideal, and a vector of integers containing the position of
the places above $P$ in the list supplied by \verb"C[3]"; if the
point lies at the infinity, the ideal is replaced by an
homogeneous irreducible polynomial in two variables.
\medskip

A \textbf{place} is represented by the four list previously cited:
a base point (\verb"list POINTS" of homogeneous coordinates); a
local equation (\verb"list LOC_EQS") for the curve at the base
point; a Hamburger-Noether expansion of the corresponding branch
(\verb"list BRANCHES"); and a local parametrization
(\verb"list PARAMETRIZATIONS") of such a branch. 
\medskip

A \textbf{divisor} is represented by a vector of integers, where
the integer at the position $i$ means the coefficient of the
$i$-th place in the divisor. 
\medskip

\textbf{Rational functions} are represented by ideals with two
homogeneous generators, the first one being the numerator of the
rational function, and the second one being the denominator.
\medskip

Furthermore, we can compute a complete list containing all the
non-singular affine (closed) places with fixed degree $d$ just by
using the procedure \verb"NSplaces" in this way:
\begin{verbatim}
> C=NSplaces(1..d,C);
\end{verbatim}

Closer to our aim is the procedure \verb"Weierstrass", which computes the non-gaps of the Weierstra\ss~semigroup at one
point and the associated functions with poles. It contains three inputs:
\begin{itemize}
    \item an \emph{integer} indicating the rational place in which we
    compute the semigroup;
    \item an \emph{integer} indicating how many non-gaps we want to
    calculate;
    \item the curve given in form of a \emph{list} \verb"C=Adj_div(f)" for
    some polynomial $f$ representing the local equation of
    the curve at the point given in the first entry.
    \end{itemize}

This procedure needs to be called from the ring \verb"C[1][2]".
Moreover, the places must be necessarily \emph{rational}.

\subsection{Procedures generalizing to several
points}\label{section:seve}

We present now a main procedure to compute the dimension of the
so-called Riemann-Roch vector spaces of the form
$\mathcal{L}(\mep)\setminus \mathcal{L}((\me-\ei)\underline{P})$.
If this dimension is equal to $1$, the procedure is also able to
compute a rational function belonging to the space.
\vspace{0.25cm}

The technique developed here is not by using the adjunction theory
directly, as we have developed theoretically in the Chapter 3
(Algorithm \ref{algo:no}), because of its high cost, but we use
the Algorithm \ref{algo:si}, or, more properly speaking, a slight
variant of it: we order the poles in a vector from the biggest one
to the smallest one (in absolute value) and we take the first in
such a vector.
{\small
\begin{verbatim}
proc RRquot (intvec m, list P, list CURVE, int chart)
"USAGE:RRquot( m, P, CURVE, ch );  m,P intvecs, CURVE a list and
ch an integer. RETURN:   an integer 0 (dimension of
L(m)\L(m-e_i)), or a list with three entries:
  @format
  RRquot[1] ideal (the associated rational function)
  RRquot[2] integer (the order of the rational function)
  RRquot[3] integer (dimension of L(m)\L(m-e_i))
  @end format
NOTE:     The procedure must be called from the ring CURVE[1][2],
          where CURVE is the output of the procedure @code{NSplaces}.
@*        P represents the coordinates of the place CURVE[3][P].
@*        Rational functions are represented by
          numerator/denominator
          in form of ideals with two homogeneous generators.
WARNING:  The place must be rational, i.e., necessarily
CURVE[3][P][1]=1. @* SEE ALSO: Adj_div, NSplaces, BrillNoether
EXAMPLE:  example RRquot; shows an example " {
  // computes a basis for the quotient of Riemann-Roch vector spaces L(m)\L(m-e_i)
  // where m=m_1 P_1 + ... + m_r P_r and m-e_i=m_1P_1+...+(m_i-1)P_i+...+m_r P_r,
  // a basis for the vector space L(m-e_i) and the orders of such functions, via
  //   Brill-Noether
  // returns 2 lists : the first consists of all the pole orders in
  //   increasing order and the second consists of the corresponding rational
  //   functions, where the last one is the basis for the quotient vector space
  // P_1,...,P_r must be RATIONAL points on the curve.

  def BS=basering;
  def SS=CURVE[5][1][1];
  intvec posinP;
  int i,dimen;
  setring SS;
  //identify the points P in the list CURVE[3]
  int nPOINTS=size(POINTS);
  for(i=1;i<=size(m);i=i+1)
   {
     posinP[i]=isPinlist(P[i],POINTS);
   }
//in case the point P is not in the list CURVE[3]
  if (posinP==0)
    {
      ERROR("The given place is not a rational place on the curve");
    }
  setring BS;
  //define the divisor containing m in the right way
  intvec D=zeroes(m,posinP,nPOINTS);
  list Places=CURVE[3];
  intvec pl=Places[posinP[chart]];
  int dP=pl[1];
  int nP=pl[2];

  //check that the points are rational
  if (dP<>1)
  {
    ERROR("The given place is not defined over the prime field");
  }
  int auxint=0;
  ideal funcion;
  funcion[1]=1;
  funcion[2]=1;
  
  // Brill-Noether algorithm
  list LmP=BrillNoether(D,CURVE);
  int lmP=size(LmP);
  if (lmP==1)
  {
    dimen=0;
    return(dimen);
  }
  list ordLmP=list();
  list sortpol=list();
    for (i=1;i<=lmP;i=i+1)
    {
      ordLmP[i]=orderRF(LmP[i],SS,nP)[1];
    }
    ordLmP=extsort(ordLmP);
    if (D[posinP[chart]] <> -ordLmP[1][1])
      {
    dimen=0;
        return(dimen);
      }
    LmP=permute_L(LmP,ordLmP[2]);
    funcion=LmP[1];
    dimen=1;
    return(list(funcion,ordLmP[1][1],dimen));
} example
 {
  "EXAMPLE:"; echo=2;
  int plevel=printlevel;
  printlevel=-1;
  ring s=2,(x,y),lp;
  poly f=y2+y3+x3+x4;
  list C=Adj_div(f);
  C=NSplaces(1,C);
  def pro_R=C[1][2];
  setring pro_R;
  intvec m=4,6;
  intvec P1=0,1,1;
  intvec P2=0,1,0;
  list P=P1,P2;
  int chart=1;
  RRquot(m,P,C,chart);
  printlevel=plevel;
 }
\end{verbatim}
}
Let us see an example:
{\small
\begin{verbatim}
> example RRquot;
// proc RRquot from lib brnoeth.lib
EXAMPLE:
  int plevel=printlevel;
  printlevel=-1;
  ring s=2,(x,y),lp;
  poly f=y2+y3+x3+x4;
  list C=Adj_div(f);
The genus of the curve is 2
  C=NSplaces(1,C);
  def pro_R=C[1][2];
  setring pro_R;
  intvec m=4,6;
  intvec P1=0,1,1;
  intvec P2=0,1,0;
  list P=P1,P2;
  int chart=1;
  RRquot(m,P,C,chart);
Vector basis successfully computed 
-->[1]:
   _[1]=x3+yz2
   _[2]=xyz+xz2
-->[2]:
   -4
-->[3]:
   1
  printlevel=plevel;
\end{verbatim}
}
This procedure needs also the following auxiliar procedures:
\vspace{0.25cm}

As \verb"RRquot" reads off the point through its homogeneous
coordinates we need to localize that point in the list
\verb"POINTS" and make the correspondence between such a point and
its position in the list of points contained in the third output
of the procedure \verb"Adj_div". This is done by mean of the
routine \verb"isPinlist". Its inputs are the point $P$ in
homogeneous coordinates, that is, a vector of integers, and the
list $L$ of points from \verb"Adj_div". The output is an integer
being zero if the point is not in the list or a positive integer
indicating the position of $P$ in $L$. Look at the example:
{\small
\begin{verbatim}
> example isPinlist;
// proc isPinlist from lib brnoeth.lib
EXAMPLE:
  ring r=0,(x,y),ls;
  intvec P=1,0,1;
  list POINTS=list(list(1,0,1),list(1,0,0));
  isPinlist( P,POINTS);
-->1
\end{verbatim}
}

We need also a procedure for ordering a list of integers. This is
partially solved by the procedure \verb"sort" from
\verb"general.lib". But \verb"sort" is not able to order lists of
negative numbers, so we have extended this algorithm to
\verb"extsort". The \verb"extsort" procedure needs to permute a
vector of integers by the instructions given by another similar
vector. This is actually done for lists of integers
(\verb"permute_L" in \verb"brnoeth.lib"), but not for vectors of
integers. This lack is covered by the procedure \verb"perm_L",
whose entries are a pair of vectors, the second vector fixing the
permutation of the first one. The output consists of the
permutated vector, as the following example shows:
{\small
\begin{verbatim}
> example extsort;
// proc extsort from lib brnoeth.lib
EXAMPLE:
  ring r=0,(x,y),ls;
  list L=10,9,8,0,7,1,-2,4,-6,3,0;
  extsort(L);
-->[1]:
   -6,-2,0,0,1,3,4,7,8,9,10
-->[2]:
   9,7,4,11,6,10,8,5,3,2,1
\end{verbatim}
}

Finally, it was interesting to fix the system for reading off the
data of the divisor needed in the \verb"BrillNoether" procedure.
Our routine zeroes takes two vectors of integers \verb"m" and
\verb"pos", and an integer \verb"siz" and it builds up a vector of
size \verb"siz", with the values contained in \verb"m" set in the
places given by \verb"pos" and zeroes in the other places. This
algorithm is called \verb"zeroes":
{\small
\begin{verbatim}
> example zeroes;
// proc zeroes from lib brnoeth.lib
EXAMPLE:
  ring r=0,(x,y),ls;
  intvec m=4,6;
  intvec pos=4,2;
  zeroes(m,pos,5);
-->0,6,0,4,0
\end{verbatim}}


\end{document}